\documentclass[a4paper]{article}

\usepackage{amsmath,amstext}  
\usepackage{amsfonts,amssymb, bbm} 
\usepackage[LGR, T1]{fontenc}
\usepackage[utf8]{inputenc}
\usepackage[english]{babel}
\usepackage{color}
\usepackage{graphicx, graphics}
\usepackage{subfigure}
\usepackage{epsfig}
\usepackage{epsf}
\usepackage{vmargin}
\newcommand{\ind}[1]{\mathbbm{1}_{\{#1\}}}
\newcommand{\wh}{\widehat}

\newcommand{\wt}{\widetilde}

\newtheorem{theo}{Theorem}[section]

\newtheorem{exam}[theo]{Example}

\usepackage[color,final]{showkeys}

\renewcommand{\c}{\textcolor{cyan}}
\newcommand{\removefornow}[1]{}

\title{\LARGE \bf  Linear minimum mean square filters for Markov jump linear systems}
\author{ Eduardo~F.~Costa and Beno\^\i te~de~Saporta
\thanks{Eduardo F. Costa is with Univ. S\~ao Paulo -
Instituto de Ci\^{e}ncias Mathem\'{a}ticas e de Computa\c{c}\~{a}o, C.P. 668,
13560-970, S\~{a}o Carlos, SP, Brazil  {\tt\small efcosta@icmc.usp.br.}
Beno\^\i te de Saporta is with University of Montpellier, 
F-34095 Montpellier, France, 
CNRS, IMAG, UMR 5149, F-34095 Montpellier, France
and Inria Bordeaux Sud Ouest, team CQFD, F-33400 Talence, France. {\tt\small e-mail: Benoite.de-Saporta@umontpellier.fr}
This work was supported by Inria associate team CDSS and ANR grant Piece, FAPESP and CNPq.}%
}

\begin{document}

\maketitle
\thispagestyle{empty}
\pagestyle{empty}

\begin{abstract}
New linear minimum mean square estimators are introduced in this paper 
by considering a cluster information structure in the filter design.
The set of filters constructed in this way can be ordered in a lattice 
according to the refines of clusters of the Markov chain, 
including the linear Markovian estimator at one end (with only one 
cluster) and the Kalman filter at the other hand (with as many clusters 
as Markov states).
The higher is the number of clusters, the heavier are pre-compuations 
and smaller is the estimation error, so that the cluster cardinality allows for a trade-off 
between performance and computational burden.
In this paper we propose the estimator, give the formulas for pre-computation of gains, 
present some properties, and give an illustrative numerical example. 
\end{abstract}

\section{Introduction}
There is a vast number of applications benefiting from the nice properties of the Kalman 
filter (KF). Among these properties, the possibility of pre-computation of 
gains \cite{Anderson79,Miller} is of much relevance for applications. 
However, in some cases pre-computation is not possible or viable due to missing \emph{a priori} 
relevant information. This is the case when using the KF to estimate the 
state of Markov jump linear systems (MJLS), since the parameters are not 
known prior to the current time instant $k$, in fact they depend on 
the Markov chain current state $\theta(k)$. 
Then, to use KF for MJLS one needs to do either online computation of the 
gains or offline pre-computation of a number of sample path dependent gains,  
a figure that grows exponentially with time. 

This drawback of KF for MJLS is one of the main motivations behind the 
emergence of other filters for this class of systems, 
see e.g. \cite{Costa11autom_filter,CostaFragosoMarques05,Fioravanti08,Goncalves10_signalproc}. 
Among them, the one that is closer to the KF in terms of structure and performance 
is the linear minimum mean square estimator (LMMSE) first introduced 
for MJLS in \cite{Costa94}. 
Indeed, as we shall see later, the LMMSE computation relies on coupled 
Riccati equations that are quite similar to the ones arising in Kalman filtering. 
Another similarity is that both are optimal in the mean square error sense, 
although under different constraints.
The main dissimilarity with KF lies in the fact that, instead of having path dependent gains,  
the LMMSE has sets of $N$ gains, $N$ being the cardinality of the Markov state space 
(which we assume finite in this paper), that are precomputed based on the system 
matrices and the initial distribution of $\theta$. During application, after observing 
$\theta(k)$ one picks the corresponding gain from the precomputed sets of gains.
In this way, one obtains the best estimate for the system state $x_k$ 
among all estimators that are \emph{linear} and \emph{Markovian}, these being 
precisely the constraints we mentioned before. 

However, there is currently no intermediary solution between the KF and the LMMSE 
in literature. In this paper we provide a ``lattice'' of filters 
bridging the KF to the LMMSE by relaxing the markovianity constraint 
and allowing clustered information of the Markov chain to be 
considered when designing the filter gains. 
By clustered information, we mean that we have a partition of the state space of the Markov chain into several classes called clusters, and we observe the trajectory of classes the chain belongs to along time.
The clusters may be chosen as they are considered as a design parameter, 
establishing a trade-off between complexity and performance 
that can be explored in practical systems aiming at the best feasible performance.
At one extreme when only one cluster is taken into account 
our filter is equivalent to the LMMSE, and at another
extreme with $N$ clusters we retrieve the Kalman filter; 
intermediary number of clusters leads to filters with variable 
performances and computational burden. 
Reasonably enough, the higher is the number of clusters the 
smaller is the attained estimation error 
and higher is the number of gains to compute.

We start with a simple, precise formulation of the optimal estimation problem 
in Section \ref{sec-problem}, 
with the estimator in the classical form of Luenberger observers. 
We then proceed in Section \ref{sec-proof} to a constructive proof that evaluates the estimation error 
and uses the completion of squares method to obtain the optimal gains. 
Some remarks on how the proposed class of estimators includes both the KF and LMMSE, and on 
the number of gains and Riccati-like equations to be precomputed, are presented.
Some variants of the studied optimization problem and how to extend optimality 
to general estimators are briefly discussed in Section \ref{sec-properties}.
We have also included a numerical example in Section \ref{sec-example} comparing the computational 
burden in terms of CPU time and the estimation error computed both 
via Monte Carlo simulation and via the proposed formula. 
The example makes clear that the performance is strongly dependent 
on the number of clusters and how the Markov states are distributed in the 
clusters.

\section{Problem formulation}\label{sec-problem}
\noindent{}Consider the MJLS
\begin{equation}\label{eq-def-system}
\begin{aligned}
x_{k+1}&= A_{\theta(k)} x_k + G_{\theta(k)} w_k
\\ y_{k}&=   L_{\theta(k)} x_k + H_{\theta(k)} w_k, \qquad k\geq 0,
\end{aligned}
\end{equation}
with initial condition $x_0\sim N(\bar x,\Psi)$  where $N(\bar x,\Psi)$ is the normal distribution with mean $\bar x$ and covariance matrix $\Psi$.
The variable $\theta(k)$ denotes the state of a Markov chain with finite state space $\{1,2,\ldots,N\}$ and initial 
distribution $\pi_0=[\text{Pr}(\theta(0)=1)\, \cdots \text{Pr}(\theta(0)=N)]$.
The noise sequence $w$ is independent from $x_0$ and the Markov chain $\theta$, $E[w_k]=0$ and $E[w_kw_k']$ is the identity matrix for all $k$.
We assume that $G_iH_i'=0$ and $H_iH_i'>0$, $1\leq i \leq N$. 
We shall consider Luenberger observers in the form\footnote{General recursive linear estimators are briefly addressed in Section \ref{rem-general-lmmse}.}
\begin{equation}\label{eq-def-filter}
\begin{aligned}
\widehat x_{k+1}&= A_{\theta(k)} \widehat x_k + M_k (y_k-L_{\theta(k)} \widehat x_k),
\end{aligned}
\end{equation}
where matrix $M_k$ is referred to as the filter gain, and the initial 
estimate is given by $\widehat x_0 = \bar x$. 
This produces an estimation error $\wt x=x-\widehat x$ satisfying
\begin{equation}\label{eq-def-filter}
\begin{aligned}
\wt x_{k+1}&= (A_{\theta(k)}-M_k L_{\theta(k)}) \wt x_k+ (G_{\theta(k)}-M_k H_{\theta(k)}) w_k,
\end{aligned}
\end{equation}
and $\wt x_0 = x_0-\bar x \sim N(0,\Psi)$. 

As for the Markov chain, we consider a partition $S_1,\ldots, S_{N_C}$ for its state space, 
and employ the variable $\rho(k)$ to indicate the partition being visited at time $k$, 
that is, 
$$\rho(k)= \sum_{m=1}^{N_C}  m\times \ind{\theta(k)\in S_m}.$$ 
We assume that the observations of $\rho$ up to time $k$ are available 
to calculate the filter gains. We also assume that the jump variable 
$\theta$ at time $k$ is available, however in the clustered information filter 
we do not take into account its past values, that is, $\theta(0),\ldots,\theta(k-1)$ 
are not taken into account when calculating the gain 
(to avoid an excessive number of branches, as explained earlier). 
Moreover, the gain should not depend on future 
information, as it has to be implemented at every time instant $k$.
Therefore, we impose that the filter matrices at time instant $k$ are in the form
\begin{equation}\label{eq-info-constraint}
\begin{aligned}
M_k &= h_k(\rho(0),\ldots,\rho(k-1),\theta(k)), \quad k\geq 0
\end{aligned}
\end{equation}
for measurable functions $h_k$. Gains satisfying this constraint are referred to as \emph{feasible gains}.
We are interested in obtaining the minimum mean square state estimation 
at a given time $s\geq 0$, leading to the optimization problem
\begin{equation}\label{eq-def-problem}
\min_{M_0,\ldots,M_s} E\{\|x_s - \widehat x_s\|^2 |\mathcal{R}_s \},
                             \qquad  \text{s.t.~} \eqref{eq-info-constraint},
\end{equation}
where we write $\mathcal{R}_s=\{\rho(0),\ldots,\rho(s),y(0),\ldots,y(s),\theta(s)\}$ 
for ease of notation. 
We refer to the filter \eqref{eq-def-filter} satisfying \eqref{eq-def-problem} as the \emph{clustered information LMMSE}, 
or CLMMSE for short.

\section{Clustered information LMMSE computation}\label{sec-proof}
Consider the MJLS in \eqref{eq-def-system} and the filter in \eqref{eq-def-filter} 
with an arbitrary sequence of feasible gains $M=\{M_k, k\geq 0\}.$
For each $k\geq 0$, $0\leq i \leq N$ and $0\leq \ell_m \leq M$, $0\leq m\leq k-1$,  
we define
\begin{equation}\label{eq-def-P}
X_{\ell_0,\ldots,\ell_{k-1},i,k}(M)=
    E(\wt x_k \wt x_k ' \ind{\rho(0)=\ell_0,\ldots,\rho(k-1)=\ell_{k-1},\theta(k)=i } ).
\end{equation}
The variable $X$ plays an important role in the derivation of 
the formula for the optimal filter because the optimal gains $M^\star=\{M_k^\star, 0\leq k\leq s\}$ 
are such that $X_{\ell_0,\ldots,\ell_{k-1},i,k}(M^\star) = Y_{\ell_0,\ldots,\ell_{k-1},i,k}$ 
where $Y$ is the solution of a Riccati-like equation, as we shall see in the next theorem.
The physical interpretation for $X$ is that, when divided by $\text{Pr}(\rho(0)=\ell_0,\ldots,\rho(k-1)=\ell_{k-1},\theta(k)=i)$
it gives the conditional error covariance matrix, e.g. it coincides with the filtering Riccati
equation of the Kalman filter when $N=N_C$ (as many clusters as Markov states).\\

\begin{exam}
We illustrate the notation introduced in \eqref{eq-def-P} on a simple example. 
Consider the MJLS \eqref{eq-def-system} with 
$A_i=G_i=1$, $1\leq i\leq 3$, $\Psi=1$, $\bar x=0$, the clusters $S_1=\{1,2\}$, $S_2=\{3\}$, 
initial distribution $\pi_0=[0.5\; 0.3 \; 0.2]$ and probability matrix
$$P = \begin{bmatrix} 
0.5 & 0.4 & 0.1 \\
 1  &   0 & 0 \\
0.5 &   0 & 0.5 
\end{bmatrix}.$$
Consider trivial gains $M_k=0$, so that $\wt x_k = x_k$ and 
$E(\wt x_{k+1} \wt x_{k+1} ') =  E(\wt x_k \wt x_k ') + 1$. 
Moreover, the estimation error and the Markov state are 
independent because all modes are identical, then 
$X_{\ell_0,\ldots,\ell_{k-1},i,k}(M)=
    E(\wt x_k \wt x_k ') \text{Pr}(\rho(0)=\ell_0,\ldots,\rho(k-1)=\ell_{k-1},\theta(k)=i ).$
For instance, for $\ell_0=1$, $i=2$ and $k=1$ we have
$$\begin{aligned}
X_{1,2,1}&=E(\wt x_1 \wt x_1 ') \text{Pr}(\rho(0)=1,\theta(1)=2 )
= 2\, \text{Pr}(\theta(1)=2|\rho(0)=1)\text{Pr}(\rho(0)=1)
\\ &= 2 \,\text{Pr}(\rho(0)=1) \Big(\sum_{j\in S_{\ell_0}} \text{Pr}(\theta(1)=2|\rho(0)=1,\theta(0)=j)
    \text{Pr}(\theta(0)=j|\rho(0)=1)\Big) 
\\ &= 2\times 0.8\times (0.4\times (0.5/0.8)+0\times (0.3/0.8)) = 0.4.  
\end{aligned}
$$ 
Similarly, $X_{1,1,1}=1.1$, $X_{1,3,1}=0.1$,
$X_{2,1,1}= 0.2$,  $X_{2,2,1}= 0$ and $X_{2,3,1}= 0.2$. Note from  
\eqref{eq-def-P} that summing $X$ in the indexes corresponding to $\ell$ and $i$ 
we obtain $E(\wt x_k \wt x_k ')$, e.g. from the above we have $E(\wt x_1 \wt x_1 ')=2$.
\end{exam}
\subsubsection*{Coupled Riccatis of the CLMMSE}
the optimal gain sequence can be (pre-)computed based on the following sets of matrices. 
Let $Y_{i,0}=\pi_i(0) \Psi $ for each $1\leq i\leq N$.  
For each $k\geq 1$, let $0\leq m\leq k-1$ and compute for each  
$0\leq i \leq N$ and  $0\leq \ell_m \leq M$,
\begin{equation}\label{eq-def-Y}
Y_{\ell_0,\ldots,\ell_{k-1},i,k}= 
\begin{cases}
0, &  p_{\ell_0,\ldots,\ell_{k-1},i,k}=0, 
\\ 
\\ \begin{aligned}
& \sum_{j\in \wt S}  p_{j i}  \big[ A_{j}Y_{\ell_0,\ldots,\ell_{k-2},j,k-1}A_{j}'
       +p_{\ell_0,\ldots,\ell_{k-2},j} G_jG_j' 
\\&  \quad \qquad + A_{j}Y_{\ell_0,\ldots,\ell_{k-2},j,k-1}L_j'\big( L_j Y_{\ell_0,\ldots,\ell_{k-2},j,k-1} L_j'  
\\&  \quad \qquad  + p_{\ell_0,\ldots,\ell_{k-2},j,k-1} H_jH_j' \big)^{-1} 
     L_j Y_{\ell_0,\ldots,\ell_{k-2},j,k-1} A_j' \big],
  \end{aligned}   &  \text{otherwise,} 
\end{cases}
\end{equation}
where we denote $$
\begin{aligned}
p_{\ell_0,\ldots,\ell_{k-1},i,k}&=\text{Pr}(\rho(0)=\ell_0,\ldots,\rho(k-1)=\ell_{k-1},\theta(k)=i),
\\ \wt S &= \{j\in S_{\ell_{k-1}}:p_{\ell_0,\ldots,\ell_{k-2},j,k-1}\not=0\}.
\end{aligned}$$

\begin{theo}\label{theo-optimal-solution}
Given the realization of $\theta(k)$, $k\leq s$, and the 
corresponding cluster observations $\rho(0),\ldots,\rho(s-1)$, 
the gains $M^\star= \{M^\star_0,\ldots,M^\star_s\}$ 
of the LMMSE can be computed for each $k\leq s$ as
\begin{equation}\label{eq-optim-gains}
M^\star _k = 
\begin{cases}
0,  
   & \text{Pr}(\rho(0),\ldots,\rho(k-1),\theta(k))=0, 
\\  \begin{aligned}
  A_{\theta(k)} &Y_{\rho(0),\ldots,\rho(k-1),\theta(k),k} L_{\theta(k)}' 
 \big(L_{\theta(k)} Y_{\rho(0),\ldots,\rho(k-1),\theta(k),k} L_{\theta(k)}'
\\&  + \text{Pr}(\rho(0),\ldots,\rho(k-1),\theta(k)) H_{\theta(k)}H_{\theta(k)}'\big)^{-1},
         \end{aligned} & \text{otherwise,}
\end{cases}
\end{equation}
where $Y$ is given in \eqref{eq-def-Y}.
Moreover, the conditional second moment of the estimation error is given, 
for each $k\leq s$, by
\begin{equation}\label{eq-computing-P} 
X_{\rho(0),\ldots,\rho(k-1),\theta(k),k}(M^\star) 
= Y_{\rho(0),\ldots,\rho(k-1),\theta(k),k},
\end{equation}
and it is optimal in the sense that, for any gain sequence $M=\{M_0,\ldots, M_s\}$,
\begin{equation}\label{eq-optimality-P}
 X_{\rho(0),\ldots,\rho(k-1),\theta(k),k}(M^\star)
\leq X_{\rho(0),\ldots,\rho(k-1),\theta(k),k}(M),   \quad 0\leq k\leq s.
\end{equation}
\end{theo}

\textbf{Proof} 
We start showing that \eqref{eq-computing-P} and \eqref{eq-optimality-P} are true
for the gains prescribed in \eqref{eq-optim-gains}. We proceed by induction in $k$. 
For the time instant $k=0$ we have that 
the initial estimate is given by $\widehat x_0 = \bar x$, yielding 
$\wt x_0 = x_0-\bar x \sim N(0,\Psi)$ (irrespectively of the filter gains), 
hence 
\begin{equation}
X_{i,0}(M^\star)=X_{i,0}(M)= E(\wt x_0 \wt x_0 ' \ind{\theta(0)=i }) 
=E(\wt x_0 \wt x_0 ') E(\ind{\theta(0)=i }) = \Psi  \pi_i(0)
=Y_{i,0}.
\end{equation}
By the induction hypothesis we assume that  
\eqref{eq-computing-P} and \eqref{eq-optimality-P} are valid for some $0\leq k<s$.
In order to complete the induction, we now consider the time instant $k+1$.
For a filter with an arbitrary sequence of easible gains, denoted by $M=\{M_0,\ldots, M_s\}$,
and a given realization $\rho(0)=\ell_0,\ldots,\rho(k)=\ell_{k},\theta(k+1)=i$ we have
\begin{equation*}
X_{\ell_0,\ldots,\ell_k,i,k+1}(M)=
    E(\wt x_{k+1} \wt x_{k+1} ' 
       \ind{\rho(0)=\ell_0,\ldots,\rho(k)=\ell_{k},\theta(k+1)=i } ).
\end{equation*}
Note that the above quantity turns out to be zero whenever 
$p_{\ell_0,\ldots,\ell_{k},i,k+1}=0$,  irrespectively of $M$, which makes
\eqref{eq-computing-P} and \eqref{eq-optimality-P} trivially true for $k+1$ in this case.
Now, in case  $p_{\ell_0,\ldots,\ell_{k},i,k+1}\not= 0$ we write 
\begin{equation}\label{eq-def-Y-v01}
\begin{aligned}
& X_{\ell_0,\ldots,\ell_k,i,k+1}(M)=
    E(\wt x_{k+1} \wt x_{k+1} ' \ind{\rho(0)=\ell_0,\ldots,\rho(k)=\ell_{k},\theta(k+1)=i } )
\\ &=  
E\Big( \big( (A_{\theta(k)}-M_k L_{\theta(k)}) \wt x_k+ (G_{\theta(k)}-M_k H_{\theta(k)}) w_k \big)
 \\& \qquad \times \big( (A_{\theta(k)}-M_k L_{\theta(k)}) \wt x_k+ (G_{\theta(k)}-M_k H_{\theta(k)}) w_k \big)'
 \ind{\rho(0)=\ell_0,\ldots,\rho(k)=\ell_{k},\theta(k+1)=i}  \Big)
\\ &=  
E\Big( [ (A_{\theta(k)}-M_k L_{\theta(k)}) \wt x_k \wt x_k' (A_{\theta(k)}-M_k L_{\theta(k)})' 
 \\& \qquad  + (G_{\theta(k)}-M_k H_{\theta(k)}) w_k w_k' (G_{\theta(k)}-M_k H_{\theta(k)})'] 
 \ind{\rho(0)=\ell_0,\ldots,\rho(k)=\ell_{k},\theta(k+1)=i } \Big)
\\ &=  
E\Big(  \sum_{j=1}^N [(A_{\theta(k)}-M_k L_{\theta(k)}) \wt x_k 
                  \wt x_k' (A_{\theta(k)}-M_k L_{\theta(k)})' 
 \\& \qquad  + (G_{\theta(k)}-M_k H_{\theta(k)}) w_k w_k' (G_{\theta(k)}-M_k H_{\theta(k)})' ] 
 \ind{\rho(0)=\ell_0,\ldots,\rho(k)=\ell_{k},\theta(k+1)=i,\theta(k)=j } \Big)
\\ &=  
E\Big(  \sum_{j\in S_{\ell_k}} [(A_j-M_k L_j) \wt x_k 
                  \wt x_k' (A_j-M_k L_j)' 
 \\& \qquad  + (G_j-M_k H_j) w_k w_k' (G_j-M_k H_j)' ] 
 \ind{\rho(0)=\ell_0,\ldots,\rho(k-1)=\ell_{k-1},\theta(k+1)=i,\theta(k)=j } \Big)
\end{aligned}
\end{equation}
where the last equality comes from the fact that $\text{Pr}(\rho(k)=\ell_{k},\theta(k)=j)=0$ whenever 
$j$ is not in the cluster $S_{\ell_k}$, and 
$\text{Pr}(\rho(k)=\ell_{k},\theta(k)=j))=\text{Pr}(\theta(k)=j)$ otherwise. 
Resuming the above calculation and writing 
$\mathcal R=\{\rho(0)=\ell_0,\ldots,\rho(k-1)=\ell_{k-1},\theta(k+1)=i,\theta(k)=j\}$ for 
ease of notation, we have:
\begin{equation}\label{eq-def-Y-v02}
\begin{aligned}
X_{\ell_0,\ldots,\ell_{k},i,k+1}(M)
&= \sum_{j\in S_{\ell_k}} \text{Pr}(\mathcal R) \big[ (A_{j}-M_k L_{j})
 E\big(\wt x_k \wt x_k'  \big| \mathcal R)  (A_{j}-M_k L_{j})' 
\\& \qquad \qquad + (G_{j}-M_k H_{j}) 
 E\big( w_k  w_k' \big| \mathcal R)  (G_{j}-M_k H_{j})' \big].
\end{aligned}
\end{equation}
From basic properties of the Markov chain we have that 
$\rho(\ell)$ and $\theta(k+1)$ are conditionally independent 
given $\theta(k)$, for any $0\leq \ell\leq k-1$. Moreover, from \eqref{eq-def-system}, 
\eqref{eq-def-filter} and \eqref{eq-info-constraint}
it can be shown that $\wt x_k$ and $\theta(k+1)$ are conditionally independent 
given $\theta(k)$, hence we may eliminate 
$\theta(k+1)=i$ from the first conditional expectation in \eqref{eq-def-Y-v02}, leading to 

\begin{equation}\label{eq-def-Y-v02b}
\begin{aligned}
X&_{\ell_0,\ldots,\ell_{k},i,k+1}(M)= \sum_{j\in S_{\ell_k}} 
    \text{Pr}(\theta(k+1)=i|\rho(0)=\ell_0,\ldots,\rho(k-1)=\ell_{k-1},\theta(k)=j)p_{\ell_0,\ldots,\ell_{k-1},j,k}
\\& \qquad \times   \big[ (A_{j}-M_k L_{j})
 E\big(\wt x_k \wt x_k' 
        \big| \rho(0)=\ell_0,\ldots,\rho(k-1)=\ell_{k-1},\theta(k)=j) 
        (A_{j}-M_k L_{j})' 
\\& \qquad + (G_{j}-M_k H_{j}) I (G_{j}-M_k H_{j})' \big]
\\& = \sum_{j\in \wt S} p_{ji}  
 \big[ (A_{j}-M_k L_{j})
        X_{\ell_0,\ldots,\ell_{k-1},j,k}(M)  (A_{j}-M_k L_{j})' 
\\& \qquad \qquad  + p_{\ell_0,\ldots,\ell_{k-1},j,k} (G_{j}-M_k H_{j}) (G_{j}-M_k H_{j})' \big],
\end{aligned}
\end{equation}
where we denote $\wt S = \{j\in S_{\ell_k}: p_{\ell_0,\ldots,\ell_{k-1},j,k}\not=0\}$.
We now turn our attention to the optimality of $M$. 
Consider a feasible gain sequence in the form 
$$\bar M=\{M_0^\star,\ldots, M_{k-1}^\star, M_k\},$$ where $M_k$ is the variable to be determined; 
since $X_{\ell_0,\ldots,\ell_{k-1},j,k}(\bar M)$  is a function of 
$M_0^\star,\ldots, M_{k-1}^\star$ only, we can use the induction hypothesis to write
\begin{align}
X_{\ell_0,\ldots,\ell_{k-1},j,k}(\bar M) = Y_{\ell_0,\ldots,\ell_{k-1},j,k}, \label{induct-aux-01}
\\  X_{\ell_0,\ldots,\ell_{k-1},j,k}(\bar M) \leq  X_{\ell_0,\ldots,\ell_{k},j,k}(M). \label{induct-aux-02}
\end{align}
Eq. \eqref{induct-aux-02} allows to write
$(A_{j}-M_k L_{j}) (X_{\ell_0,\ldots,\ell_{k-1},j,k}(\bar M) - 
        X_{\ell_0,\ldots,\ell_{k-1},j,k}(M))  (A_{j}-M_k L_{j})'\leq 0$, 
irrespectively of $M_k$, and using \eqref{eq-def-Y-v02b} we evaluate
\begin{equation}\label{eq-def-Y-v02d}
\begin{aligned}
&X_{\ell_0,\ldots,\ell_{k},i,k+1}(\bar M) -  X_{\ell_0,\ldots,\ell_{k},i,k+1}(M)
\\&
=\sum_{j\in \wt S } p_{ji}  
 \big[ (A_{j}-M_k L_{j}) (X_{\ell_0,\ldots,\ell_{k-1},j,k}(\bar M) - 
        X_{\ell_0,\ldots,\ell_{k-1},j,k}(M))  (A_{j}-M_k L_{j})' \leq 0
\end{aligned}
\end{equation}
Also, by plugging \eqref{induct-aux-01} into \eqref{eq-def-Y-v02b}, 
\begin{equation}\label{eq-def-Y-v02c}
\begin{aligned}
X_{\ell_0,\ldots,\ell_{k},i,k+1}&(\bar M)
= \sum_{j\in \wt S} p_{ji}  
 \big[ (A_{j}-M_k L_{j})
        Y_{\ell_0,\ldots,\ell_{k-1},j,k}(\bar M)  (A_{j}-M_k L_{j})' 
\\& \qquad + p_{\ell_0,\ldots,\ell_{k-1},j,k} (G_{j}-M_k H_{j}) I (G_{j}-M_k H_{j})' \big],
\end{aligned}
\end{equation}
and, by completing squares and denoting $\Phi = L_j Y_{\ell_0,\ldots,\ell_{k-1},j,k} L_j'  
 + \text{Pr}(\rho(0)=\ell_0,\ldots,\rho(k-1)=\ell_{k-1},\theta(k)=j) H_jH_j'$ for brevity, 
we obtain
\begin{equation}\label{eq-def-Y-v03}
\begin{aligned}
& X_{\ell_0,\ldots,\ell_{k},i,k+1}(\bar M)= \sum_{j\in \wt S} p_{ji }  
 \big[ A_{j}Y_{\ell_0,\ldots,\ell_{k-1},j,k}A_{j}'
       +p_{\ell_0,\ldots,\ell_{k-1},j,k} G_jG_j' 
\\& \qquad \qquad\qquad \qquad     + \big( M_k - A_{j}Y_{\ell_0,\ldots,\ell_{k-1},j,k}L_j' \Phi^{-1} \big)
        \Phi \big( M_k - A_{j}Y_{\ell_0,\ldots,\ell_{k-1},i,k}L_j' \Phi^{-1} \big)'
\\& \qquad \qquad\qquad \qquad  - A_{j}Y_{\ell_0,\ldots,\ell_{k-1},j,k}L_j' \Phi^{-1} L_j Y_{\ell_0,\ldots,\ell_{k-1},j,k} A_j' \big],
\end{aligned}
\end{equation}
thus making clear that the minimal $X$ is attained by setting 
\begin{equation*}
\begin{aligned}
M_k&=M_k^\star=g(\ell_0,\ldots,\ell_{k-1},j) = A_{j}Y_{\ell_0,\ldots,\ell_{k-1},i,k}L_j'\Phi^{-1} 
         \\& =A_{j}Y_{\ell_0,\ldots,\ell_{k-1},i,k}L_j' \big(L_j Y_{\ell_0,\ldots,\ell_{k-1},j,k} L_i'  
              + p_{\ell_0,\ldots,\ell_{k-1},j,k} H_jH_j'\big)^{-1},
\end{aligned}
\end{equation*}
whenever $p_{\ell_0,\ldots,\ell_{k-1},j,k} \not = 0$, 
confirming the second equation in \eqref{eq-optim-gains};  
the inverse always exists because we have assumed $H_iH_i'>0$. 
If $j$ is such that $p_{\ell_0,\ldots,\ell_{k-1},j,k}= 0$ 
then the gain $M_k$ is immaterial for the error covariance, indeed 
we see from  \eqref{eq-def-Y-v02c} that such gain is not accounted for, 
so that one can pick $M_k=0$, confirming the first equation in \eqref{eq-optim-gains}.
Chosing the gain as above we get the gain sequence $M^\star= \{M^\star_0,\ldots,M^\star_k\}$ 
and 
$$X_{\ell_0,\ldots,\ell_{k},j,k+1}(M^\star) \leq X_{\ell_0,\ldots,\ell_{k},j,k+1}(\bar M),$$ 
so that \eqref{induct-aux-02} produces
$$X_{\ell_0,\ldots,\ell_{k},j,k+1}(M^\star) \leq X_{\ell_0,\ldots,\ell_{k},j,k+1}(M),$$
which confirms  \eqref{eq-optimality-P} for the time instant $k+1$. 
Substituting $M_k=M_k^\star$ in \eqref{eq-def-Y-v03} we get after some manipulations
\begin{equation*}
\begin{aligned}
X_{\ell_0,\ldots,\ell_{k},i,k+1}(M^\star)
&= \sum_{j\in \wt S} p_{ji } 
 \big[ A_{j}Y_{\ell_0,\ldots,\ell_{k-1},j,k}A_{j}'
       +p_{\ell_0,\ldots,\ell_{k-1},j,k} G_jG_j' 
\\& \quad \qquad  + A_{j}Y_{\ell_0,\ldots,\ell_{k-1},j,k}L_j' 
   \big( L_j Y_{\ell_0,\ldots,\ell_{k-1},i,k} L_j'  
\\& \quad \qquad    + p_{\ell_0,\ldots,\ell_{k-1},j,k} 
     H_jH_j' \big) L_j Y_{\ell_0,\ldots,\ell_{k-1},j,k} A_j' \big]= Y_{\ell_0,\ldots,\ell_{k},i,k+1},
\end{aligned}
\end{equation*}
which confirms \eqref{eq-computing-P} for $k+1$, thus completing the induction. 
It remains only to show the optimality of $M^\star$ in terms of \eqref{eq-def-problem}. 
This follows directly from \eqref{eq-optimality-P}, in fact, 
\begin{equation*}
\begin{aligned}
E\{(\wt x_k^{\star})'\wt x_k^\star \}
&=\sum  E((\wt x_k^{\star})' \wt x_k^\star \ind{\rho(0)=\ell_0,\ldots,\rho(k-1)=\ell_{k-1},\theta(k)=i } )
\\& =\sum  \text{tr}(X_{\ell_0,\ldots,\ell_{k-1},i,k}(M^\star))
\\& \leq \sum \text{tr}(X_{\ell_0,\ldots,\ell_{k-1},i,k}(M))
\\& =\sum   E(\wt x_k' \wt x_k  \ind{\rho(0)=\ell_0,\ldots,\rho(k-1)=\ell_{k-1},\theta(k)=i } )
=E\{ \wt x_k^{'}\wt x_k\},
\end{aligned}
\end{equation*}
where all sums are in the indexes ${0\leq i\leq N, 0\leq \ell_m \leq M, 0\leq m\leq k-1}$
and we denote by $\wt x_k^\star$ the estimation error associated with the gain $M^\star$.
\hfill $\Box$

\section{Properties of the CLMMSE}\label{sec-properties}
\subsection{Number of matrices to be computed and stored}
\label{rem-number-matrices}
For each $0\leq k \leq s-1$, we compute $NN_C^k$ matrices on the left hand side of \eqref{eq-def-Y}, 
hence we have (up to) this number of recursive Riccati equations to solve.
We also have the computation and storage of an equal number of gains. 
Then, to obtain the state estimate at time $s$, we have to store 
a total of $N (N_C^s-1)(N_C-1)^{-1}$ gains when $N_C\not =1$, and $sN$ gains otherwise.
Regarding the number of matrix inverses, one may invert each $Y$ given by 
\eqref{eq-def-Y} and store it at time step $k$ for the forthcoming iterates, 
hence we have a total of (up to) $N (M_C^{(s-1)}-1)(M_C-1)^{-1}$ inverses. 
\subsection{Filtering in the entire interval $0\leq k\leq s$}
\label{rem-variants-problem}
Note from \eqref{eq-optim-gains} that, given a realization of the Markov chain $\theta(k)$, $k\geq 0$, 
the time instant $s$ involved in the problem formulation \eqref{eq-def-problem} 
affects only the cardinality of the optimal gain sequence $M^\star$. 
More precisely, if $\{M_k^\star, 0\leq k\leq s\}$ is the gain sequence 
attaining \eqref{eq-def-problem}, and, if we replace $s$ with $\ell\leq s$ in   
\eqref{eq-def-problem} and obtain the new optimal gain sequence $\{\wt M_k, 0\leq k\leq \ell\}$ 
(considering the same Markov chain realization), then we have that 
$M_k^\star = \wt M_k$, $0\leq k \leq \min(\ell,s)$. 
This is consistent with the sense of optimality in \eqref{eq-optimality-P}, and is  
in perfect harmony with the theory of both Kalman filter 
and the standard LMMSE. 
As a consequence, the provided clustered information LMMSE is also a solution for 
the multiobjective problem
\begin{equation*}
\min_{M_0,\ldots,M_s} \{ E\{\|x_0 - \widehat x_0\|^2 |\mathcal{R}_0 \}, 
\ldots, E\{\|x_s - \widehat x_s\|^2 |\mathcal{R}_s \},
\end{equation*}
or for any linear combination of mean square errors writen in the form
$$\min_{M_0,\ldots,M_s} \left(\sum_{0\leq k\leq s} \alpha_k E\{\|x_k - \widehat x_k\|^2 |\mathcal{R}_k \}\right).$$
\subsection{Linking the Kalman filter and the standard LMMSE}
\label{rem-variants-problem}
It is simple to see that we retrieve the standard LMMSE  
when we consider only one partition $S_1=\{1,\ldots,N\}$. In fact, 
in this setup we have $p_{\ell_0,\ldots,\ell_{k-1},i,k}=P(\theta(k)=i)$
and one can check by inspection that \eqref{eq-def-Y} and the 
LMMSE Riccati equation \cite[Eq. XX]{Costa94} are identical.  
As for the Kalman filter, if we set $S_i=\{i\}$, $1\leq i\leq N$, then 
$$\begin{aligned}
p_{\ell_0,\ldots,\ell_{k-1},i,k}&=\text{Pr}(\rho(0)=\ell_0,\ldots,\rho(k-1)=\ell_{k-1},\theta(k)=i)
\\ &=\text{Pr}(\theta(0)=\ell_0,\ldots,\theta(k-1)=\ell_{k-1},\theta(k)=i) 
= \pi_{\ell_0}(0)p_{\ell_0,\ell_1}\cdots p_{\ell_{k-1},i},
\\ \wt S &= \{j\in S_{\ell_{k-1}}:p_{\ell_0,\ldots,\ell_{k-2},j,k-1}\not=0\}= \{\ell_{k-1}\}, 
\end{aligned}$$
if $p_{\ell_0,\ldots,\ell_{k-2},\ell_{k-1},k-1}\not=0$
and \eqref{eq-def-Y} is reduced to
\begin{equation*}
\begin{aligned}
Y_{\ell_0,\ldots,\ell_{k-1},i,k}= 
&  \big[ A_{j}Y_{\ell_0,\ldots,\ell_{k-2},j,k-1}A_{j}'
       +p_{\ell_0,\ldots,\ell_{k-2},j} G_jG_j' 
\\&  \quad \qquad + A_{j}Y_{\ell_0,\ldots,\ell_{k-2},j,k-1}L_j'\big( L_j Y_{\ell_0,\ldots,\ell_{k-2},j,k-1} L_j'  
\\&  \quad \qquad  + p_{\ell_0,\ldots,\ell_{k-2},j,k-1} H_jH_j' \big)^{-1} 
     L_j Y_{\ell_0,\ldots,\ell_{k-2},j,k-1} A_j' \big].
  \end{aligned}   
\end{equation*}
From Theorem \ref{theo-optimal-solution} we have 
\begin{equation*}
\begin{aligned}
Y_{\ell_0,\ldots,\ell_{k-1},i,k}&=X_{\ell_0,\ldots,\ell_{k-1},i,k}(M^\star)
\\=&
    E(\wt x_k^\star \wt x_k^{\star '} |\theta(0)=\ell_0,\ldots,\theta(k-1)=\ell_{k-1},\theta(k)=i ) 
    \pi_{\ell_0}(0) p_{\ell_0,\ell_1}\cdots p_{\ell_{k-1},i}
\end{aligned}
\end{equation*}
so that, writing $Z_{\ell_0,\ldots,\ell_{k-1},i}=E(\wt x^\star_k \wt x\star_k' |\theta(0)=\ell_0,\ldots,\theta(k-1)=\ell_{k-1},\theta(k)=i )$, substituting in the above equation for $Y$ and manipulating (cancelling the $p$s and $\pi$s)
yields
\begin{equation}
\begin{aligned}
Z_{\ell_0,\ldots,\ell_{k-1},i,k}= 
&  \big[ A_{j}Z_{\ell_0,\ldots,\ell_{k-2},j,k-1}A_{j}'
       +G_jG_j' 
\\&  \quad \qquad + A_{j}Z_{\ell_0,\ldots,\ell_{k-2},j,k-1}L_j'\big( L_j Z_{\ell_0,\ldots,\ell_{k-2},j,k-1} L_j'  
\\&  \quad \qquad  + H_jH_j' \big)^{-1} 
     L_j Z_{\ell_0,\ldots,\ell_{k-2},j,k-1} A_j' \big],
  \end{aligned}   
\end{equation}
which is the usual Riccati difference equation appearing in Kalman filters. 
This means that $Y_{\ell_0,\ldots,\ell_{k-1},i,k}$
is equal to the Kalman covariance matrix multiplied by the probability 
that the Markov chain visits $\ell_0,\ldots,\ell_{k-1},i$. 
In the gain formula \eqref{eq-optim-gains}, this probability is cancelled, 
yielding that the Kalman gain coincide with $M_k^\star$. Concluding, we have 
the Kalman filter and the markovian LMMSE in opposite ``extremes'' of the CLMMSE, and a 
lattice of estimators between them, depending on how the Markov states are arranged in 
clusters.
\subsection{General LMMSE}
\label{rem-general-lmmse}
Consider linear estimators of the general form
\begin{equation}\label{eq-def-linfil}
\begin{aligned}
z_{k+1}&= F_{k} z_k + \bar G_{k}y_k,
\end{aligned}
\end{equation}
where matrices $F_k$ and $\bar B_{k}$, $k\geq 0$, are the optimization variables 
replacing $M_k$ in the problem \eqref{eq-def-problem}; 
consider also that $F_k = f_k(\rho(0),\ldots,\rho(k-1),\theta(k))$
and $\bar G_k = g_k(\rho(0),\ldots,\rho(k-1),\theta(k))$, where $f_k,g_k$ are measurable functions.
It can be demonstrated that the optimal estimate satisfies 
$$z_k=\wh x^\star_k, \quad 0\leq k\leq s, \;\text{a.s.},$$
which is produced by setting $\bar G_k = M_k^\star$ and 
$$F_k=A_{\theta(k)}-M_k^\star L_{\theta(k)}$$ 
where $M_k^\star$ is as in \eqref{eq-optim-gains}, 
thus retrieving the Luenberger observer form \eqref{eq-def-filter} and the solution given in Theorem \ref{theo-optimal-solution}.
This is not surprising since the ``innovation form'' $F_k=A_{\theta(k)}-M_k L_{\theta(k)}$ for some $M_k$
is necessary for some basic properties of an observer to be fulfilled, 
e.g. for $\widetilde x_k=z_k-x_k$ to remain zero a.s. for $k\geq 1$ in cases when $z_0=x_0$ a.s. and 
there is no additive noise in the state ($G_i=0$, $i=1,2,\ldots,N$).

\section{Illustrative example}\label{sec-example} 

\begin{figure}[ht]
\begin{center}
\includegraphics[width=10cm]{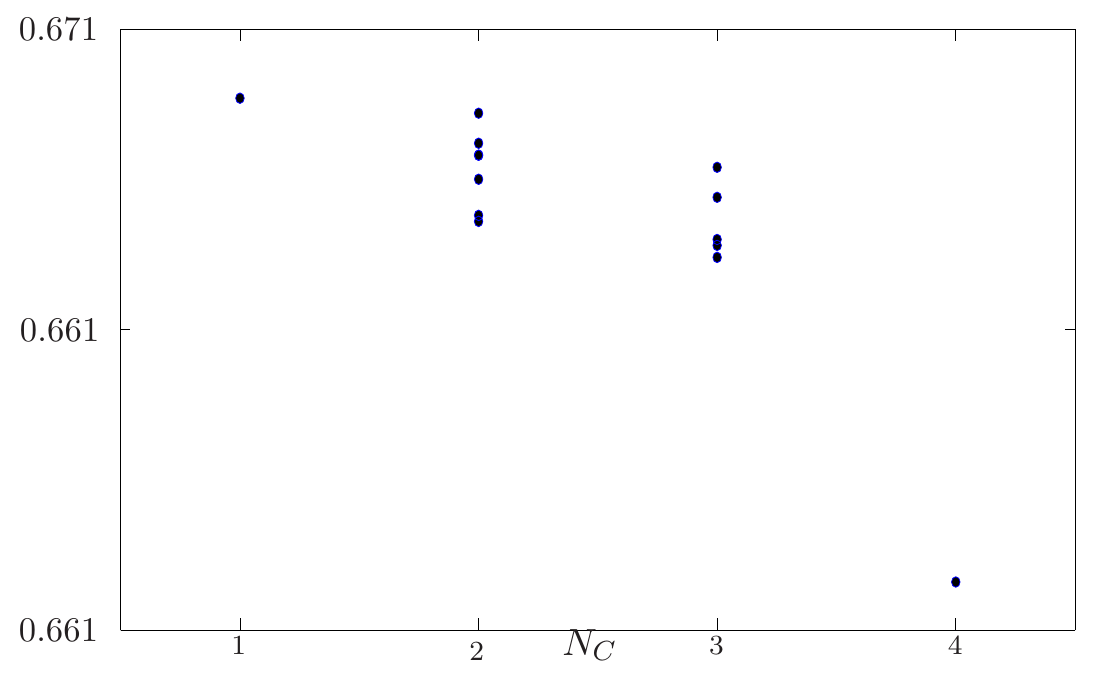}
\caption{Mean square error $E(\|\tilde x_{10}\|^2)$ versus number of clusters $N_C$ for every cluster configuration for the system
with data given in (\ref{data1}).}
\label{fig1}
\end{center}
\end{figure}

We have applied the CLMMSE to the system given in \cite{ZB09}. The system data is reproduced below for ease of reference.
\begin{equation*}
A_1=\left[\begin{array}{ll}0 -0.405\\0.81 0.81\end{array}\right],
A_2=\left[\begin{array}{ll}0 -0.2673\\0.81 1.134\end{array}\right],
A_3=\left[\begin{array}{ll}0 -0.81\\0.81 0.972\end{array}\right],
A_4=\left[\begin{array}{ll}0 -0.1863\\0.81 0.891\end{array}\right],
\end{equation*}
\begin{equation*}
G_i=\left[\begin{array}{ll}0.5&0\\0&0\end{array}\right],
L_i=\left[\begin{array}{ll}1&0\end{array}\right],
H_i=1, i=1,\ldots,4,
\end{equation*}
\begin{equation}\label{data1}
P=\left[\begin{array}{llll}
0.3& 0.2& 0.1& 0.4\\
0.3& 0.2& 0.3& 0.2\\
0.1&0.1&0.5&0.3\\
0.2&0.2&0.1&0.5
\end{array}\right].
\end{equation}
Every possible cluster configuration has been taken into account, and the aim was to estimate the state at time instant $s=10$. The mean square error was calculated by means of 
(\ref{eq-def-P}) and (\ref{eq-computing-P}), which lead to  $E(\|\tilde x_{10}\|^2)=\sum_{\ell_0,\ldots,\ell_9,i}\text{trace}(Y_{\ell_0,\ldots,\ell_9,i,10})$; results have been confirmed by Monte Carlo simulation. Figure~\ref{fig1} shows the obtained results in groups according to the number of clusters $N_C$. As expected, the standard LMMSE (with $N_C=1$) presented the largest estimation error and the Kalman filter ($N_C=4$) features the smallest one. The performance of other filters with "intermediary configurations" ($N_C=2,3$) is similar to the LMMSE regarding the error, hence they are not much appealing in view of their higher complexity. Note that in this particular example, the modes are similar to each other (only two parameters of $A$ change).

\begin{figure}[ht]
\begin{center}
\includegraphics[width=10cm]{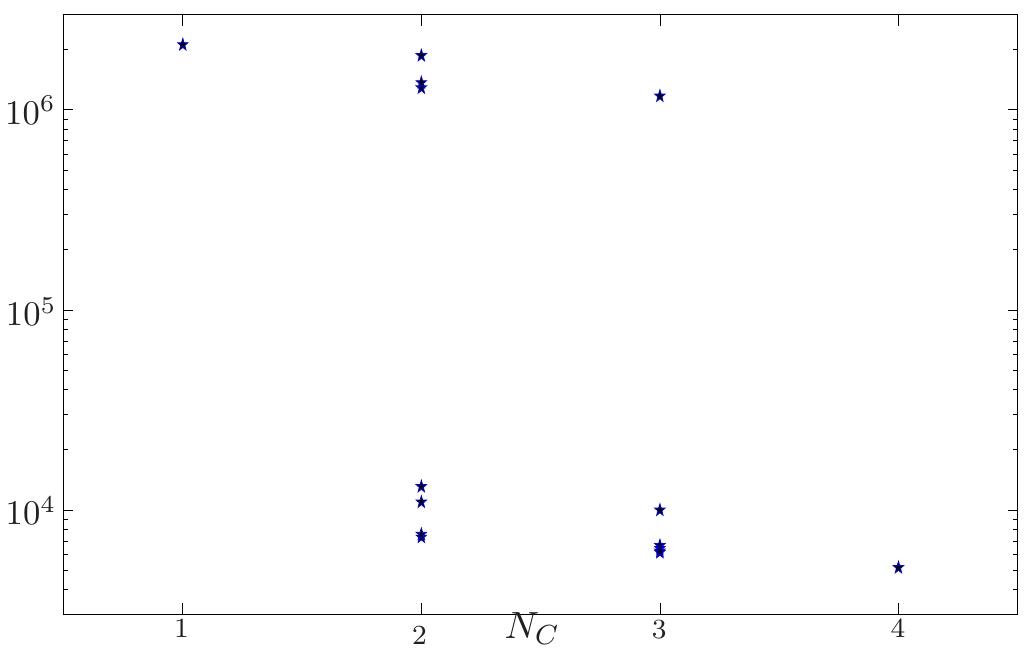}
\caption{Mean square error $E(\|\tilde x_{10}\|^2)$ versus number of clusters $N_C$ for every cluster configuration for the system
with data given in (\ref{toto}).}
\label{fig2}
\end{center}
\end{figure}

Now let us introduce a more relevant change in one mode by replacing $A_4$, $G_4$, $L_4$ and $H_4$ with
\begin{equation}\label{toto}
10A_4, 10G_4, 10L_4 \text{ and } 10H_4,
\end{equation}
respectively. The results are displayed in Figure~\ref{fig2}. We now can clearly distinguish two groups of
filters, one with average errors around 106 and a second one around 104. In this setup, there is a
tendency for better performance when $\theta=4$ is isolated from other states, e.g. with $S_1 = \{1, 2, 3\}$
and $S_2 = \{4\}$ we have $E(\|\tilde x_{10}\|^2)\simeq 7, 566.69$, while $S_1 = \{1, 2\}$ and $S_2 = \{3, 4\}$ lead to
$E(\|\tilde x_{10}\|^2)\simeq1, 286, 306.18$. Moreover, considering that the (hard to implement) Kalman filter
yields $E(\|\tilde x_{10}\|^2)\simeq 5, 152.95$ we see that the filter with configuration $S_1 = \{1, 2, 3\}$
and $S_2 = \{4\}$ is quite competitive in this scenario.

\section{Concluding remarks}
We have explored the Markov state information structure in MJLS leading to new filters
whose estimation error and complexity lie in the between the standard Kalman filter and the
LMMSE, and establish a trade-off between performance and computational burden. This allows
one to explore the computational resources at hands more deeply, seeking for the best possible
estimates. As illustrated in the example, the new filters can provide competitive alternatives to
the existing ones. We note that, for a given plant, the computational burden depends only on
$N_C$, see the formulas given in Section~\ref{rem-number-matrices}. Then, a relevant question is how to select the
Markov states to form each cluster so that the estimation error is minimized, which will be
considered in future research.


\end{document}